\documentclass[a4paper,10pt]{article}
\title{On curvature and feedback classification of two-dimensional optimal control systems}
\author{Ulysse Serres\footnote{Universit\'e de Bourgogne UFR des Sciences et des Techniques
9 avenue Alain Savary, BP 47870-21078 Dijon cedex France and
SISSA/ISAS Via Beirut 2-4, 34014 Trieste Italy; email: serres@sissa.it}}
\date{}

\usepackage{amsmath}
\usepackage{amsthm}
\usepackage{amssymb}

\newtheorem{theorem}{Theorem}[section]
\newtheorem{definition}[theorem]{Definition}
\newtheorem{example}[theorem]{Example}

\renewcommand{\le}{\leqslant}
\renewcommand{\ge}{\geqslant}

\newcommand{\lang}{\left\langle}
\newcommand{\rang}{\right\rangle}
\newcommand{\bs}{\boldsymbol}
\newcommand{\hamvec}[1]{\boldsymbol{\vec{#1}}}
\newcommand{\parfrac}[2]{\frac{\partial{#1}}{\partial{#2}}}
\newcommand{\R}{\mathbb{R}}

\numberwithin{equation}{section}

\begin{document}

\maketitle

\begin{abstract}
The goal of this paper is to extend to two-dimensional optimal control systems with scalar input the classical
notion of Gaussian curvature of two-dimensional Riemannian surface using the Cartan's moving frame method.
This notion was already introduced by A. A. Agrachev and R. V. Gamkrelidze for more general control systems
using a purely variational approach.
Then we will see that the ``control'' analogue to Gaussian curvature reflects similar intrinsic properties
of the extremal flow. In particular if the curvature is negative, arbitrarily long segment of extremals are
locally optimal. Finally, we will define and characterize flat control systems.
\end{abstract}

\section{Introduction}
In Riemannian geometry the Gaussian curvature of a manifold reflects intrinsic properties of the geodesic flow,
i.e. properties that do not depend on the choice of local coordinates.
For example, the geodesics of the surface have no conjugate points if the curvature is non-positive.
Indeed, these geodesics are extremals of a particular time optimal control problem the dynamics of which is given
by
\begin{equation*}
\dot{q}=\cos u\,\bs{e}_1(q)+\sin u\,\bs{e}_2(q),\quad u\in S^1,
\end{equation*}
where $(\bs{e}_1,\bs{e}_2)$ forms an orthonormal frame of the Riemannian structure on the manifold.
Our goal is to generalize the classical notion of Gaussian curvature of two-dimensional Riemannian surfaces
for two-dimensional smooth optimal control problems.
The notion of curvature tensor for non linear optimal control problems was first introduced in \cite{AAA1}
by A. A. Agrachev and R. V. Gamkrelidze with a purely variational description by means of Jacobi curves, which
are curves in the Lagrangian Grassmannian.
Here we will not deal with Jacobi curves but use the Cartan's moving frame method in order to construct a
feedback invariant frame associated to our optimal control problem and provide a less general but also more
geometric definition of the curvature function.

Consider a control system of the form
\begin{equation}\label{cs}
\dot{q}=\bs{f}(q,u),\quad q\in M,\quad u\in U,
\end{equation}
where $M$ and $U$ are smooth connected manifolds. Let $\bs{\tilde{f}}(\tilde{q},\tilde{u})$, 
$(\tilde{q},\tilde{u})\in \tilde{M}\times \tilde{U}$,
be the right-hand side of another such system. We say that the two systems
are {\it feedback-equivalent} if there exists a diffeomorphism $\Theta:M\times U\to\tilde{M}\times\tilde{U}$ of
the form
\begin{equation}\label{feedback}
\Theta(q,u)=(\phi(q),\psi(q,u))
\end{equation}
which transforms the first system to the second, i.e. such that
\begin{equation*}
T_q\phi(\bs{f}(q,u))=\bs{\tilde{f}}(\phi(q),\psi(q,u)).
\end{equation*}
In the above diffeomorphism $\phi$ plays the role of a change of coordinates in the state space $M$, and $\psi$
called {\it pure feedback transformation} reparametrizes the set of controls $U$ in a way depending on the state
variable $q\in M$.
Our aim is to provide feedback invariants for control system (\ref{cs}) when the manifold $M$ is of dimension 
two and the control set $U$ of dimension one what we suppose from now.

In this case, if the coordinates on the manifold are fixed, a control system of type (\ref{cs}) is parametrized
by two functions of three variables, and the group of feedback transformations of type (\ref{feedback}) is
parametrized by two functions of two variables and one function of three variables. Therefore, we can a priori
normalize only one function among the two functions defining control system (\ref{cs}). Thus, we expect to have 
only $2-1=1$ ``principal" feedback invariant, i.e. a function of three variables, in this equivalence problem.

All results of the present paper will be presented without proof. Anyway, most of these proofs can be found in 
the references cited at the end.

\section{Curvature}
Suppose that we want to minimize an integral cost
$\int_{t_0}^{t_1}\varphi(q,u)dt$, along the trajectories of control
system
(\ref{cs}). We write the normal maximized Hamiltonian function of PMP (Pontryagin Maximum Principle) which is
defined by
\begin{equation}\label{hamPMP}
h(\lambda)=\max_{u \in U}\left(\lang \lambda,\bs{f}(q,u) \rang-\varphi(q,u)\right),
\quad \lambda \in T^*_q M,\quad q\in M,
\end{equation}
where $\lang \cdot\;\!,\cdot \rang$ denotes the canonical pairing between the tangent and the cotangent bundles
over $M$. Hamiltonian $h$ is a function on the cotangent bundle $T^*M$ which is, because of being independent
of $u$, feedback-invariant. Thus, all objects construct from Hamiltonian $h$ through intrinsic relations will
also be feedback invariants.
As usual, if $\sigma$ denotes the symplectic two-form of the cotangent bundle over $M$, we define, via the
relation $\sigma(\cdot,\hamvec{h})=dh$, the Hamiltonian vector field $\hamvec{h}$ associated to the 
Hamiltonian function $h$.
Assume that $h$ is a smooth function, then the corresponding Hamiltonian vector field $\hamvec{h}$
is well-defined and tangent to the level set of $h$.
PMP asserts (see e.g. \cite{AAAbook}) that optimal trajectories of system (\ref{cs}) are projections onto $M$ of
trajectories of the Hamiltonian system $\dot{\lambda}=\hamvec{h}(\lambda)$, in other words trajectories of 
Hamiltonian field $\hamvec{h}$ are extremals of our optimal control problem.
Now fix a level set $\mathcal{H}=h^{-1}(e)$ of our Hamiltonian, then the intersection
$\mathcal{H}_q=\mathcal{H}\cap T^*_qM$ is a curve in the plane $T^*_qM$ and under the regularity assumptions
\begin{equation}\label{reg}
\bs{f}(q,u) \wedge \parfrac{\bs{f}(q,u)}{u} \neq 0,\quad
\parfrac{\bs{f}(q,u)}{u} \wedge \parfrac{^2\bs{f}(q,u)}{u^2}>0,
\quad q \in M,\quad u\in U,
\end{equation}
such a curve is a strictly convex curve surrounding the origin and it admits, up to sign and translation, 
a natural parameter providing us with a vector field $\bs{v}_q$ on $\mathcal{H}_q$ and by consequence with a 
vertical vector field $\bs{v}$ on $\mathcal{H}$.
Vector field $\bs{v}$ is characterized by the fact that it is, up to sign, the unique vector field on
$\mathcal{H}$ such that
\begin{equation}\label{b}
L^2_{\bs{v}}s = -s + b L_{\bs{v}}s,
\end{equation}
where $s$ denotes the restriction to $\mathcal{H}$ of Liouville one-form ``$pdq$" of $T^*M$ and $b$ is a smooth
function on the level $\mathcal{H}$.

Actually function $b$ is the feedback invariant of our control system that characterizes Riemannian problems.
Namely control problem (\ref{cs}) defines a Riemannian  geodesic problem if and only if the invariant $b$ is 
identically equal to zero.

Since vector fields $\hamvec{h}$ and $\bs{v}$ are feedback-invariant, it is natural to think that the curvature
of our system may arise from a commutator relation of these fields. Indeed, the following theorem confirms this
intuition.
\begin{theorem}\label{curvatura}
Vector fields $\bs{v}$ and $\hamvec{h}$ satisfy the following nontrivial commutator relation:
\begin{equation}\label{courbure}
\Big[ \hamvec{h},\Big[\bs{v},\hamvec{h}\Big]\Big]=\kappa \bs{v}.
\end{equation}
\end{theorem}
Proof of this theorem can be found in \cite{AAAbook}, \cite{ulysse}.
The coefficient $\kappa$ in the identity (\ref{courbure}) is defined to be the {\it curvature} of our optimal 
control problem and since the fields $\hamvec{h}$ and $\bs{v}$ are feedback-invariant, the curvature
$\kappa$ is also feedback-invariant.

Fix a system of local coordinates $\lambda=(\theta,q)\in\mathcal{H}$ where $\theta$ parametrizes the fiber
$\mathcal{H}_q$ so that $\bs{v}=\parfrac{}{\theta}$ and denote $\parfrac{}{\theta}=\,'\,$.
Define the function $c=c(\theta,q)$ by 
\begin{equation}\label{FunctionC}
d_qs=c s\wedge s',
\end{equation}
where $d_qs$ is the differential of
the Liouville one-form $s$ with respect to the horizontal coordinates. Then, the Hamiltonian field takes the form
\begin{equation*}
\hamvec{h}=\bs{f}-c\parfrac{}{\theta},
\end{equation*}
and the curvature $\kappa$ is evaluated as follows:
\begin{equation}\label{evcourboure}
\kappa(\theta,q)=L_{\hamvec{h}'}c-L_{\hamvec{h}}c'.
\end{equation}
\begin{example}{\rm
Consider the control system corresponding to the geodesic problem on a two-dimensional Riemannian manifold:
\begin{equation}
\dot{q}=\cos u\,\bs{e}_1(q)+\sin u\,\bs{e}_2(q),\quad u\in S^1.
\end{equation}
In this case, control curvature $\kappa$ is the {\it Gaussian curvature} of the Riemannian manifold $M$ and it is
evaluated as follows:
\begin{equation}\label{riemcurv}
\kappa(q)=-c_1^2 -c_2^2 + L_{\bs{e}_1}c_2 - L_{\bs{e}_2}c_1,
\end{equation}
where $c_1$, $c_2$ are the structural constants of the orthonormal frame $(\bs{e}_1,\bs{e}_2)$ on $M$:
\begin{equation*}
[\bs{e}_1,\bs{e}_2]=c_1\bs{e}_1+c_2\bs{e}_2,\quad c_1,\ c_2\in C^\infty(M).
\end{equation*}
See \cite{AAAbook} for the proof of this formula. Of course, for the Riemannian problem the curvature
$\kappa=\kappa(q)$ depends only on the base point $q\in M$ as one can see from formula (\ref{riemcurv}) but in
general this is not the case: the curvature $\kappa$ depends also on the coordinate in the fiber $\mathcal{H}_q$ 
and thus is a function on the whole three-dimensional manifold $\mathcal{H}$.
}\end{example}

Observe that relations (\ref{b}) and (\ref{courbure}) define two feedback invariants: the function $b$ and the
curvature $\kappa$. Both $b$ and $\kappa$ are functions on the three-dimensional level surface $\mathcal{H}$, so
that they are principal feedback invariants of our control system. Since our feedback equivalence problem admits
only one invariant these functions are not ``independent". Indeed invariants $b$ and $\kappa$ are connected by
the following differential relation:
\begin{equation}\label{bnk}
L_{\bs{v}}\kappa+b\kappa+L^2_{\hamvec{h}}b=0,
\end{equation}
which can easily be derived calculating some bracket relations between vector fields $\bs{v}$ and $\hamvec{h}$.
In particular, relation (\ref{bnk}) shows that in the special case of Riemannian problems, the curvature
$\kappa$ is a function on the base manifold $M$ without any computation. Indeed, since Riemannian problems
are characterized by the vanishing of function $b$, (\ref{bnk}) reduces to $L_{\bs{v}}\kappa=0$.

\section{Jacobi equation}
It is easy to see that the regularity assumptions (\ref{reg}) imply
$\hamvec{h}\wedge\bs{v}\wedge[\hamvec{h},\bs{v}]\neq 0$
so that vector fields $\hamvec{h}$, $\bs{v}$, $[\hamvec{h},\bs{v}]$ form a moving frame on the level surface
$\mathcal{H}$.
In this section we use this moving frame to derive an ODE on conjugate time of our two-dimensional optimal 
control problem. This ODE, Jacobi equation in the moving frame, will show that the control curvature analogue 
to the Gaussian curvature enjoys similar properties.

Fix a point $q_0\in M$ and define a two-dimensional surface in $\mathcal{H}$ by:
\begin{equation*}
\mathcal{L}^t_0=e^{t\,\hamvec{h}}(\mathcal{H}_{q_0}),\quad t\in \R,
\end{equation*}
where $e^{t\,\hamvec{h}}$ denotes the flow of the Hamiltonian field $\hamvec{h}$. The surface $\mathcal{L}^t_0$
is the lift in the cotangent bundle of trajectories $t\mapsto q(t)$ in $M$ of control system (\ref{cs}) with 
starting point $q(0)=q_0$. We say that a point $q=q(t)$, $t\neq 0$, is conjugate to $q_0$ 
(or time $t$ is conjugate to zero) if $q$ is a critical value of the canonical projection
\begin{equation}\label{conjproj}
\pi:\mathcal{L}^t_0 \to M.
\end{equation}
It is easy to check that the tangent space $T_\lambda \mathcal{L}_0^t$, $\lambda\in\mathcal{L}^t_0$,
is spanned by the vectors $\hamvec{h}(\lambda)$ and $(e^{t\,\hamvec{h}}_*\bs{v})(\lambda)$ so that the point 
$q(t)=\pi(\lambda)$ is conjugate to $q_0$ if and only if
\begin{equation*}
(e^{t\,\hamvec{h}}_*\bs{v})(\lambda)\in {\rm span}\left(\hamvec{h}(\lambda),\bs{v}(\lambda)\right).
\end{equation*}
Consider the decomposition of the vector field $e^{t\,\hamvec{h}}_*\bs{v}$ in our moving frame on $\mathcal{H}$:
\begin{equation*}
e^{t\,\hamvec{h}}_*\bs{v}=\alpha(t)\hamvec{h}+\beta(t)\bs{v}+\gamma(t)\Big[\hamvec{h},\bs{v}\Big].
\end{equation*}
It turns out that coefficients $\alpha(t)$, $\beta(t)$, $\gamma(t)$ are solutions to the Cauchy problem
\begin{equation}\label{cauchy}
\left(\begin{array}{c}
      \dot{\alpha} \\
      \dot{\beta} \\
      \dot{\gamma}
      \end{array}
\right)=
\left(\begin{array}{ccc}
      0 & 0 & 0 \\
      0 & 0 & \kappa_t \\
      0 & -1  & 0
      \end{array}
\right)
\left(\begin{array}{c}
      \alpha \\
      \beta \\
      \gamma
      \end{array}
\right),\quad \alpha(0)=1,\quad \beta(0)=\gamma(0)=0,
\end{equation}
where $\kappa_t=\kappa(e^{t\,\hamvec{h}}(\lambda_0))$, $\pi(\lambda_0)=q_0$. It is quite obvious that Cauchy 
problem (\ref{cauchy}) is equivalent to the second order linear ODE, called {\it Jacobi equation}
\begin{equation}\label{boundary}
\ddot{\gamma}+\kappa_t\gamma=0,\quad \gamma(0)=\gamma(t)=0,
\end{equation}
where $\kappa_t=\kappa(e^{t\,\hamvec{h}}(\lambda_0))$.
Thus an instant $t$ is a conjugate time for our optimal
control problem if and only if there exists a non trivial solution to the boundary value problem (\ref{boundary}).
Using the Sturm's comparison theorem for second order ODEs one can prove the following theorem about the 
occurrence of conjugate points for system (\ref{cs}).
\begin{theorem}\label{comparison}
Let $q(t)$, $q(0)=q_0$, be a solution of an optimal two-dimensional control problem and let $\kappa_t$
be the value of the curvature along an extremal $\lambda(t)$, $\pi(\lambda(t))=q(t)$.
\renewcommand{\theenumi}{\roman{enumi}}
\renewcommand{\labelenumi}{$(\theenumi)$}
\begin{enumerate}
\item If $\kappa_t\le 0$ for all $t\ge 0$, then $q_0$ has no conjugate points for $t\in[0,+\infty]$.
\item If $\kappa_t\le \kappa_1$ (resp. $\kappa_t<\kappa_1$) for all $t\ge 0$, and some constant $\kappa_1>0$,
then $q_0$ has no conjugate points along $q(\cdot)$ for $t\in[0,\pi/{\sqrt{\kappa_1}\,}[$
(resp. for $t\in[0,\pi/{\sqrt{\kappa_1}\,}]$).
\item If $0<\kappa_0\le \kappa_t$ (resp. $0<\kappa_0< \kappa_t$), for all $t\ge 0$,
then $q_0$ must have at least a conjugate point for $t\in\;]0,\pi/{\sqrt{\kappa_0}\,}]$
(resp. for $t\in\;]0,\pi/{\sqrt{\kappa_0}\,}[\,$).
\end{enumerate}
\end{theorem}
The following theorem gives sufficient condition for a trajectory $q(t)$ on $M$ to be strongly locally optimal 
in terms of conjugate points (see \cite{AAAbook} for the definition of strong optimality and the proof of the
following theorem).
\begin{theorem}
Let the trajectory $q(t)$ be as in theorem \ref{comparison}. If the time interval $]0,t_1]$ does not contain
conjugate points, then the trajectory $q(t)$ is strongly locally optimal for $t\in [0,t_1]$.
\end{theorem}
On the other hand if an instant $t_c \in\;]0,t_1]$ is conjugate to zero, then there exists an instant
$\tilde{t}\in\;]0,t_1]$ where the trajectory $q(t)$, $t\in\;]0,t_1]$, ceases to be locally optimal.
\begin{example}\label{zermelo}
{\rm Zermelo navigation problem (see \cite{cara} for a detailed description).
This problem is a time optimal control problem which consists of finding the
quickest nautical path of a yacht in the presence of stationary sea currents. The sea surface is modeled by a
two-dimensional Riemannian surface $M$ and the currents by an autonomous vector field $\bs{X}\in {\rm Vec\,}M$.
Dynamics of optimal trajectories for Zermelo problem are given by
\begin{equation*}
\dot{q}=\bs{X}(q)+\cos u\,\bs{e}_1+\sin u\,\bs{e}_2,\quad q\in M,\quad u\in S^1,
\end{equation*}
where $(\bs{e}_1,\bs{e}_2)$ form an orthonormal frame of the Riemannian surface $M$.
Suppose that the manifold $M$ is the Euclidean plane $\R^2$. Then, in the coordinate system $(q_1,q_2,u)$
vector fields $\hamvec{h}$ and $\bs{v}$ read
\begin{eqnarray*}
\hamvec{h}&=&(\bs{X}_1+\cos u)\parfrac{}{q_1}+(\bs{X}_2+\sin u)\parfrac{}{q_2}
-\big\langle D_q\bs{X}\textstyle{{{-\sin u}\choose{\phantom{-}\cos u}}},
\textstyle{{{\cos u}\choose{\sin u}}}\big\rangle \displaystyle{\parfrac{}{u}},\\
\bs{v}&=&\sqrt{\bs{X}_1(q)\cos u+\bs{X}_2(q)\sin u+1}\parfrac{}{u},
\end{eqnarray*}
so that one can compute the curvature using formula (\ref{courbure}) (here, $\lang \cdot\;\!,\cdot \rang$ denotes 
the scalar product between vectors). If we suppose moreover that the drift term is the linear field
$\bs{X}(q)={{\ \,a\ b}\choose{-b\ a}}q$ then, the control curvature for this problem is $\kappa=-a^2/4$
(see \cite{ulysse}) and theorem \ref{comparison} thus implies that there is no conjugate point along trajectories.
The following theorem can also be proved using the definition of conjugate points.
\begin{theorem}
There is no conjugate point for Zermelo navigation problem on $\R^2$
when the drift term $\bs{X}$ is a linear vector field.
\end{theorem}
Let us sketch the proof of this result
(see \cite{ulysse} for the detailed proof).
Since the drift term is linear, it is easy to compute the map $\pi$ (see (\ref{conjproj}) for the definition)
which takes the form:
\begin{equation*}
\pi(t,u_0)=q(t,q_0,u_0)=
e^{tA}q_0+\int_{0}^{t}e^{(t-\tau)A}\textstyle{ {{\cos u(\tau)}\choose{ \sin u(\tau)}} }d\tau,
\end{equation*}
where $A$ is the matrix representation of the linear drift term.
Now, because the Hamiltonian flow preserves the Liouville one-form, it is easy to see that the
differential $d_{(t,u_0)}\pi$ is of maximal rank if and only if $d_{u_0}\pi$ 
is of maximal rank. Saying this, it is now an easy task to find a vector $v$
such that $d_{u_0}\pi(v)\neq 0$ which completes the proof.
The above theorem is also valid on $\R^n$ where the Zermelo navigation problem can be generalized without any
difficulty; the proof is also similar.

}\end{example}

\section{Flat systems}
In Riemannian geometry it is well-known that if the Gaussian curvature of the surface is nonzero then, one can
not rectify simultaneously the geodesics by a change of coordinates. Only Riemannian flat systems, i.e. systems for
which the geodesics are ``straight lines" have this property. For control systems the situation is quite 
different first of all because control systems with zero curvature are not necessarily flat. We present here
a new theorem which gives a characterization of flat control systems in terms of the feedback invariants
$\kappa$ and $b$. We begin with the following definition.
\begin{definition}
A control system $\dot{q}=\bs{f}(q,u)$ is said to be flat if it is feedback equivalent to a control system of
the form $\dot{q}=\bs{f}(u)$.
\end{definition}
It is obvious that a flat system has zero curvature but the contrary is in general not true. For example a 
Zermelo problem defined on  the Euclidean plane $\R^2$ with a nonzero linear drift term is never flat.

Suppose that a control system satisfies 
\begin{equation}\label{Lhb=0}
L_{\hamvec{h}}b=0.
\end{equation}
The above property implies in particular that the plane curves $\mathcal{H}_q\subset T^*M$ are all
of the same centro-affine length.
Control systems of this type are very peculiar and have nice geometric properties that we do not discuss here.
However such systems with zero curvature are characterized in the theorem below.
\begin{theorem}\label{classif1}
There exists a feedback transformation such that:
\begin{equation}\label{bracket}
\left[\bs{f}(\cdot,u),\parfrac{\bs{f}(\cdot,u)}{u} \right]=0
\end{equation}
if and only if the feedback invariants $\kappa$ and $L_{\hamvec{h}}b$ are identically equal to zero.
Moreover if we fix local coordinates $q=(q_1,q_2)$ in $M$, then these systems can be parametrized by a
one-parameter family of diffeomorphisms generated by the vector field:
\begin{equation}\label{Xu}
\bs{X}_u=(a_1(u)+q_2)\parfrac{}{q_1}+(a_2(u,q_2)-q_1)\parfrac{}{q_2}.
\end{equation}
\end{theorem}
In the above theorem if $u$ is a control parameter such that the fields $\bs{f}$ and $\parfrac{\bs{f}}{u}$ commute
then, vector field $\bs{X}_u$ is the infinitesimal generator of a diffeomorphism
$P_u\in{\rm Diff\,}M$ such that
\begin{equation}\label{Frob}
P_{u*}\left(\bs{f}(\cdot,u),\parfrac{\bs{f}(\cdot,u)}{u} \right)=
\left(\left(\begin{array}{c} 1 \\ 0 \end{array}\right),
      \left(\begin{array}{c} 0 \\ 1 \end{array}\right)\right).
\end{equation}
Notice that commutativity between vector fields $\bs{f}$ and $\parfrac{\bs{f}}{u}$ is not a feedback-invariant
property.
When the curvature is identically zero the above theorem shows that the PDE (\ref{Lhb=0}) can be reduced to the
nonautononous ODE
\begin{equation*}
\frac{dq}{du}=\bs{X}_u(q).
\end{equation*}
The following theorem characterizes flat control systems.
\begin{theorem}\label{classif2}
A control system of type (\ref{cs}) is flat if and only if its feedback invariants
$\kappa$, $L_{\hamvec{h}}b$ and $L_{[\boldsymbol{v},\hamvec{h}]}b$ vanish identically.
\end{theorem}
We do not discuss in detail proofs of theorems \ref{classif1} and \ref{classif2} in this paper but we
roughly explain the main ideas.
The proofs are based on the following differential equation which can easily be derived from the
differentiation of the structural equations of our  feedback invariant moving frame on $\mathcal{H}$:
\begin{equation}\label{PDEforC}
c''+bc'+c=L_{\hamvec{h}}b,
\end{equation}
where $c$ is the function defined in (\ref{FunctionC}).
It follows immediately from this equation that if a control system is such that (\ref{bracket}) holds
(respectively if a control system is flat) then, its feedback invariants $\kappa$ and $L_{\hamvec{h}}b$
(respectively $\kappa$, $L_{\hamvec{h}}b$ and $L_{[\bs{v},\hamvec{h}]}b$) vanish identically.
To prove the converse observe first that if a control system has zero curvature then,
the vector fields $\hamvec{h}$ and $[\bs{v},\hamvec{h}]$ commute so that the choice of a natural parameter
$\theta$ on the fibers $\mathcal{H}_q$ defines a foliation of $\mathcal{H}$, the leaves of which are formed
by the trajectories of the fields $\hamvec{h}$ and $[\bs{v},\hamvec{h}]$.
Now, choose the parameter $\theta$ (recall that this natural parameter is fixed only
up to transformation of the form $\theta\mapsto\pm\theta+\phi(q)$) so that $c$ becomes zero which is possible 
since $c$ satisfies equation (\ref{PDEforC}) with $L_{\hamvec{h}}b=0$.
This shows in particular that there exists a feedback transformation so that (\ref{bracket}) holds and by the
Frobenius theorem one gets the existence of a diffeomorphism $P_u\in {\rm Diff\,}M$ such
that (\ref{Frob}) holds.
In order to get the expression (\ref{Xu})
we use Moser's argument for which the key idea is to determine the diffeomorphisms $P_u$ by representing them as
the flow of a family of vector fields $\bs{X}_u$ on $M$. We thus suppose that
\begin{equation*}
\frac{d}{dt}P_u=\bs{X}_u\circ P_u,\quad P_{u_0}={\rm Id},
\end{equation*}
and the expression of $\bs{X}_u$ in coordinates follows from differentiation with respect to $u$ of
(\ref{Frob}). This complete the proof of theorem \ref{classif1}. In order to complete the proof of theorem
\ref{classif2} one has just to check that $L_{\hamvec{h}}b=0$ and $L_{[\bs{v},\hamvec{h}]}b=0$ imply that
$b=b(u)$ which, in addition with (\ref{bracket}) and $\kappa=0$ easily implies that the system is flat.

We now conclude our discussion with the following example.
\begin{example}{\rm
Consider Zermelo navigation problem as in example \ref{zermelo}. One can prove that this problem is flat if
and only if the Riemannian surface $M$ is flat and the drift term $\bs{X}$ is constant.
}\end{example}

\subsection*{Acknowledgments}
I am grateful to Professor Andrei A. Agrachev for fruitful discussions.


\begin{thebibliography}{99}
\bibitem{AAA1}
A. A. Agrachev, R. V. Gamkrelidze,
Feedback-Invariant Optimal Control Theory and Differential Geometry-I. Regular Extremals,
Journal of Dynamical and Control Systems, {\bf 3} (1997), 343-389.
%
\bibitem{AAAbook}
A. A. Agrachev, Yu. L. Sachkov,
{\it Control Theory from the Geometric Viewpoint},
Springer-Verlag, 2004.
%
\bibitem{cara}
C. Carath\'eodory,
{\it Calculus of Variations},
Chelsea Publishing Company, 1989 (third edition), New York, {\bf \S 276-\S 460}.
%
\bibitem{ulysse}
U. Serres
{\it On the curvature of two-dimensional optimal control systems and Zermelo's navigation problem},
to appear in Journal of Mathematical Sciences.
\end{thebibliography}
\end{document}